\documentclass[10pt]{article}
\usepackage{amssymb,ulem,amsmath,cancel,mathrsfs, capt-of}
\usepackage[colorlinks=true, pdfstartview=FitV, linkcolor=blue, citecolor=blue, urlcolor=blue]{hyperref}

\usepackage{geometry} 
\usepackage{framed,color}
\definecolor{shadecolor}{rgb}{0.9, 0.9, 0.81}
\usepackage{ tikz, pgfplots, todonotes,ulem,amssymb}
\usetikzlibrary{decorations.markings,arrows, calc, shapes, intersections, patterns}

\def \scr{\mathscr}

\def \ds{\displaystyle}
\def\le{\left}

\def\ri{\right}

\def\bea#1\eea{\begin{align}#1\end{align}}

\def \&{\hspace{-15pt}&}
\def \d{{\mathrm d}}

\newtheorem{definition}{Definition}[section]
\def \bd{\begin{definition}}
\def \ed{\end{definition}}
\def \bp#1\ep{
 \definecolor{shadecolor}{rgb}{0.95, 0.95, 0.86}
 \begin{shaded}\begin{proposition} #1 \end{proposition}\end{shaded}}

\def \bt#1\et{
 \definecolor{shadecolor}{rgb}{0.95, 0.95, 0.86}
 \begin{shaded}\begin{theorem} #1 \end{theorem}\end{shaded}}

\def \bl#1\el{
 \definecolor{shadecolor}{rgb}{0.95, 0.95, 0.86}
 \begin{shaded}\begin{lemma} #1 \end{lemma}\end{shaded}}

\newtheorem{problem}{Problem}[section]
\newtheorem{theorem}{Theorem}[section]
\newtheorem{proposition}{Proposition}[section]

\newtheorem{remark}{Remark}[section]
\newtheorem{lemma}{Lemma}[section]

\def\K{\mathcal K}

\def\nn{\nonumber}
\def\be{\begin{equation}}
\def\ee{\end{equation}}
\def\ben{\begin{displaymath}}
\def\een{\end{displaymath}}

\def\baa{\begin{eqnarray}}
\def\eaa{\end{eqnarray}}

\def\ba{\begin{array}}
\def\ea{\end{array}}
\makeatletter
\@addtoreset{equation}{section}
\makeatother

\def \eqref #1{(\ref{#1})}
\def \1{\mathbf 1}
\def \br{\begin{remark}}
\def\er{\end{remark}}

\def\C{{\mathbb C}}
\def\Z{{\mathbb  Z}}
\def\R{{\mathbb R}}

\def\Q{{\mathcal Q}}
\def\N{{\mathbb N}}

\def\ov{\overline}
\def \Re{\mathrm {Re}\,}
\def \Im {\mathrm {Im}\,}
\def \E{\mathcal E}

\begin{document}

\title{Numerical computation algorithm}

\vspace{0.2cm}
\begin{center}
\begin{Large}
\bf Chebotarov continua, Jenkins-Strebel differentials and related problems: a numerical approach
\end{Large}

\bigskip
M. Bertola$^{\dagger}$\footnote{Marco.Bertola@concordia.ca},  
\bigskip
\begin{small}

$^{\dagger}$ {\it   Department of Mathematics and
Statistics, Concordia University\\ 1455 de Maisonneuve W., Montr\'eal, Qu\'ebec,
Canada H3G 1M8} \\
\end{small}
\vspace{0.5cm}
\end{center}

%

\begin{abstract}
We detail a numerical algorithm and related code to construct rational quadratic differentials on the Riemann sphere that satisfy the Boutroux condition. These differentials, in special cases, provide solutions of (generalized) Chebotarov problem as well as being instances of Jenkins--Strebel differentials. The algorithm allows to construct Boutroux differentials with prescribed polar part, thus being useful in the theory of weighted capacity and Random Matrices. 

\end{abstract}
\tableofcontents

\section{Introduction}
The classical Chebotarov problem consists in the following: 
\begin{problem}
\label{chebo}
Given $\E=\{e_1,\dots, e_N\}\subset \C$,  consider the class consisting  of  continua (connected compact set) $\K$  such that $\E\subset \K$. 
The problem is to find a set $\K_0$ minimizing the capacity  ${\rm Cap}(\K)$ within said class, and showing its uniqueness.
\end{problem}
The problem was suggested by P\'olya to Chebotarov in a letter,   and was solved by Lavrentieff \cite{Laurentiev2}.
Goluzin then proved \cite{Goluzin} that the extremal continuum is the collection of all critical horizontal trajectories of the rational quadratic differential  on the Riemann sphere of the form 
\be
\label{Goluzin}
Q(z)\d z^2 = -\frac {\prod_{\ell=1}^{N-2} (z-d_\ell)}{\prod_{\ell=1}^N(z-e_\ell)} \d z^2.
\ee
We refer to the monograph of Strebel \cite{StrebelBook} for all the definitions and properties of trajectories of quadratic differentials. 
We only recall that a path $\gamma$ is a horizontal trajectory of $Q(z)\d z^2$ in the neighbourhood of a regular point (where $Q(z)\neq 0$) if $ \Im\le( \sqrt{Q(\gamma(t))} \dot\gamma(t)\ri)=0$: note that the condition depends  neither on the parametrization of $\gamma$ nor on the choice of determination of the square root. A {\it critical} (horizontal) trajectory is one such that the closure of its maximal extension contains a zero of arbitrary multiplicity  or a simple pole of $Q$. 
  The points $d_1,\dots, d_{N-2}$  need to be determined implicitly by the problem. Goluzin  also characterized in terms of the uniformizing map 
\be
\zeta: \{|w|>1\} \to \C\setminus \E,
\ee
showing that $\zeta(w)$ must satisfy the differential equation 
\be
(w \zeta'(w))^2 = \frac {\prod_{\ell=1}^N(\zeta(w)-e_\ell)}{\prod_{\ell=1}^{N-2} (\zeta(w)-d_\ell)}.
\ee
See also a review in \cite{Ortega}, where the authors propose an approximation scheme based on the above theorem.

The connection with the theory of quadratic differential is thus  well known in the community of approximation and geometric-function theorists.  Very recently the generalized Chebotarov problem (with external field) has been addressed in \cite{GuilhermeAlvez}, based on the previous work \cite{GuilhermeKuijlaars}. 

A very relevant notion, which is not generally well appreciated in the above communities, is that of {\it Jenkins--Strebel (quadratic) differential}. We refer to \cite{StrebelBook} and the review \cite{ArbarelloCornalba}.  The definition requires some preparation; let $\scr R$ be a Riemann surface of genus $g\geq 0$, $\E = \{e_1,\dots, e_N\}\subset \scr R$ a collection of distinct points and let $\mathcal P = \{p_1,\dots, p_K\}$ another distinct collection. To each $p_j$ we associated a real positive number $a_j>0$ (the ``perimeter'').  We denote $\dot {\scr R}:= \scr R\setminus \mathcal P$ the pointed surface.

The notion of Jenkins--Strebel differential is then the following definition.
\bd[\cite{ArbarelloCornalba}, Def 31]
\label{defJS}
The quadratic differential $\Q$  on $\dot{\scr R }$ is a Jenkins--Strebel differential if 
\begin{enumerate}
\item It has double poles at each $p_j$ with an expression in terms of any local coordinate near $p_j$of the form 
\be
\Q = -\frac{a_j^2}{4\pi^2}\frac {\d z^2}{z^2} (1+ \mathcal O(z));
\ee
\item $\Q$ extends to $\scr R $ to a differential with  at most simple poles at $\E$;
\item all the non critical horizontal trajectories are closed.
\end{enumerate}
\ed
 The theorem of Jenkins and Strebel below is  that such a quadratic differential exists and is unique for any choices of $\E,  P$ and perimeters $a_j>0$: 
\begin{theorem}[\cite{ArbarelloCornalba} Theroem 32]
For any choices of sets $\mathcal P, \E$ and perimeters $a_j>0$ there is a unique Jenkins--Strebel differential as in Def. \ref{defJS} such that, in addition, the complement of the critical graph $\Gamma$ is a union of $K$ disjoint disks $\mathbb D_j$ with $p_j\in\mathbb D_j$.
\end{theorem}

The relationship between the Chebotarov problem and the Jenkins--Strebel differential is the following (simply established) theorem, which is however not well recognized in the literature.
\begin{theorem}
The quadratic differential solving the Chebotarov problem is the Jenkins--Strebel differential for the data $\scr R  = \mathbb P^1$, $K=1$, $p_1=\infty$, $a_1 = 2\pi$ and $\E  = \{e_1,\dots, e_N\}$. 
\end{theorem}
Indeed, in this case the complement of the critical graph $\Gamma$ is a single conformal disk containing $\infty$, and $\E\subset \Gamma$, with $\Gamma$ itself being the extremal set. See also \cite{MFRakh}, \cite{JenkinsArt}.

The application of Jenkins--Strebel-like quadratic differential is also pervasive in the asymptotic analysis of Riemann--Hilbert problems associated to Random Matrix Theory since one key ingredient for the technique is the construction of a ``$g$--function'', which is, in most instances, nothing but a special case of antiderivative of a Boutroux differential.

The paper presents an algorithm, and the accompanying Matlab code, to find these so--called Boutroux differentials. 
In particular these differentials address precisely the critical measures discussed in \cite{GuilhermeAlvez}.

\paragraph{Boutroux quadratic differentials.}
The notion was proposed in  \cite{BertoBoutroux} to address the construction of the $g$--function in the asymptotic analysis of non-hermitean orthogonal polynomials.  We state the definition in the case of a meromorphic quadratic differential on the Riemann sphere $\mathbb P^1$, but the definition is easily generalized.
\bd
\label{boutroux}
A meromorphic quadratic differential $\Q(z) = Q(z) \d z^2$ on $\mathbb P^1$, with $Q(z)$ a rational function,  is a {\bf Boutroux} differential if the differential $\eta(z) = \sqrt{\Q(z)} = \sqrt{Q(z)} \d z$ on the (hyperelliptic) Riemann surface $\scr R$ of the square root of $Q(z)$ has all periods purely real; 
\be
\oint_\gamma \eta \in \R
\ee
where $\gamma$ is any closed contour on $\scr R$ that avoids the poles of $\eta$.
\ed
The horizontal trajectories are the level-sets of the imaginary part of an antiderivative of $\eta$, $V(\zeta) = \Im \int^\zeta \eta$, which is single--valued and harmonic away from the poles of $\eta$. Note that the simple poles of $Q$ are regular points of $\eta$ on $\scr R$. Double poles of $\Q$ become simple poles of $\eta$ and the Boutroux condition requires that the residue of $\eta$ at these points is imaginary. 

Since all Jenkins--Strebel differentials are also Boutroux differentials, our focus now shifts on the latter. 
The Boutroux condition is meaningful also in the presence of poles of $\Q$  of order higher than $2$; in this case there are horizontal trajectories that do not close but end at one of the higher order poles of $\Q$. Nonetheless the horizontal trajectories determine a regular {\it foliation} in view of the fact that the function $V(\zeta)$ is  a well defined and single--valued harmonic function on $\scr R$. 

 \paragraph{Main goal.}
The goal of the paper is  simply stated as follows 
\be
\text{Determine a numerically implementable algorithm to find Boutroux differentials.} 
\ee
We are not going to seek the maximal generality and we will content ourselves in showing how to find quadratic differentials with an arbitrary number of simple poles in the finite region of the plane and a pole of even order at $z=\infty$, with prescribed singular part. We trust that the skillful reader could adapt the ideas explained in the next section to even more general situations. 

\paragraph{Acknowledgements.}
The work of the author was supported in part by the Natural Sciences and Engineering Research Council of Canada (NSERC) grant RGPIN-2023-04747.

\section{Boutroux differentials with prescribed polar part at $\infty$}
\label{algorithm}
The goal of the section is to show how to numerically solve the following problem 
\begin{problem}
\label{problemPhiE}
Let  $\mathcal E=\{e_1,\dots, e_N\} \subset \C$ be a collection of distinct points, $R, L\in \N$  and $\Phi(z) = \sum_{\ell=1}^R \frac{t_\ell}\ell z^\ell$ a polynomial with complex coefficients. Let $t_0\in \R$ arbitrary.  Determine a rational quadratic differential $\Q(z) = -Q(z)\d z^2$ of the form 
\be
\label{formQ}
Q(z) = \frac {S^2(z) \Delta(z)} {E(z)}, \ \ \ \deg S = L, \ \  \deg \Delta = M, \ \  2L+M - N = 2(R-1),
\ee
where $E(z)$ is the polynomial with roots at $\mathcal E$, the polynomial $S$ has leading coefficient $t_R$,  the polynomial  $\Delta(z)$ is monic,   $S,\Delta$ are generically prime, $\Delta$ has simple zeros and 
\begin{enumerate}
\item one determination of the square root of $Q$ has the behaviour
\be
\sqrt{Q(z)} \d z=\le( \Phi'(z) + \frac {t_0}z + \mathcal O(z^{-2}) \ri)\d z , \ \ \ |z|\to \infty ;\label{singpart}
\ee
\item on the Riemann surface $\scr R$ determined by $\eta^2 = Q(z)$, of genus (generically) $g = (M+N)/2-1 = R+L+N-2$, we have the {\bf Boutroux condition}\footnote{This is the same Boutroux condition as in Def. \ref{boutroux} because the differential $\Q$ is $-Q(z)\d z$. The extra minus implies that the periods  of $\sqrt Q \d z$ are now imaginary rather than real. We found this formulation less confusing.}
\be
\oint_\gamma \sqrt{Q} \d z \in i\R, \ \ \ \forall \gamma\in H_1(\scr R\setminus \{\infty_\pm\}, \Z),
\ee
namely, all periods of the differential $\sqrt{Q} \d z$ on $\scr R$ minus the two points above $z=\infty$ are imaginary. 
\end{enumerate}
\end{problem}
The setup is that of \cite{GuilhermeAlvez} for a simple case where the only higher order pole is at infinity. 
We now describe an algorithm of gradient descent that can be implemented numerically to obtain instances of these Boutroux differentials.

The main idea is to begin by  choosing  a random $Q(z)$ of the form \eqref{formQ} and deform it along the gradient lines of the positive functional
\be
\label{defF}
\scr F[Q]:=\frac 1 2 \sum_{\ell=0}^R \le|t_\ell-\oint_{|z|>>1}\frac{\sqrt {Q(z)} \d z}{2i\pi z^\ell} \ri|^2 +\frac 1 2 \sum_{\ell=1}^{2g} \Re\le(\oint_{\mathcal \gamma_\ell} \sqrt{Q} \d z\ri)^2.
\ee
Observe that $\scr F[Q]$ is zero exactly if and only if $Q$ solves the Problem \ref{problemPhiE}.

\paragraph{General abstract considerations.}
For a given fixed $\E$ the functional $\scr F$ can be considered as a positive function on a moduli space consisting of  the Riemann surface of $\sqrt Q$ parametrized by the roots of $S$, $\Delta$ (or their coefficients) and by the {\it marking} of the surface, namely, the choice of a complete basis of homology cycles $\gamma_1,\dots, \gamma_{2g}$. Note that if we change the marking in general the value of $\scr F$ will change; however if we are at a zero of $\scr F$ where all the real parts of the contour integrals along $\gamma_j$ are zero, then this property will be independent of the choice of marking.

This moduli space is smooth if we remove the ``discriminant locus'' defined as the locus where roots of $\Delta$ are not simple or coincide with some of the roots of $E$ or $S$. The functional $\scr F$ has no critical point other than the points of absolute minimum where $\scr F =0$ in the interior of this moduli space. This follows from the determinacy of the linear system \eqref{sysf0} below, from which it follows that the gradient of  $\scr F$ never vanishes. 

It can (and does) happen that the minimum is achieved at some point along the boundary, for example where a root of $\Delta$ coincides with a point of $\E$. In this case we have a ``drop in genus'' since the limiting hyperelliptic surface is of lower genus than the one we started. Ditto if two roots of $\Delta$ coalesce. These are "non-generic" situations in the space of parameters $t_\ell$ and configurations $\E$; in the random-matrix literature these correspond to transition regimes and are actually of significant interest. However, for the purpose of this perfunctory investigation we content ourselves to work under genericity assumptions. 

\begin{figure}
\begin{minipage}{0.48\textwidth} 
\includegraphics[width=1\textwidth]{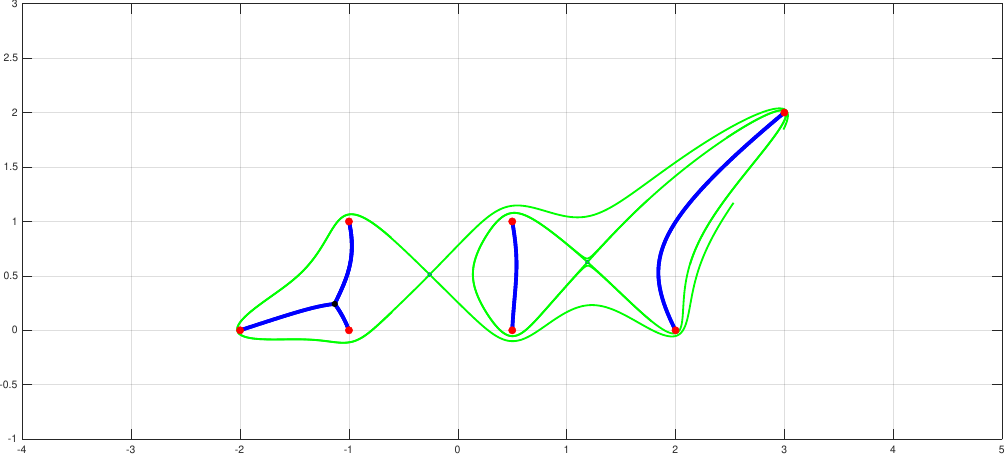} 
{\small Chebotarov continuum for the set $E=\{-1+0i,0.5+0i,0.5+1i,-1+1i,2+0i,-2+0i,3+2i\}$ with $L=2$.} 
\end{minipage} 
\hfill
\begin{minipage}{0.48\textwidth} 
\includegraphics[width=1\textwidth]{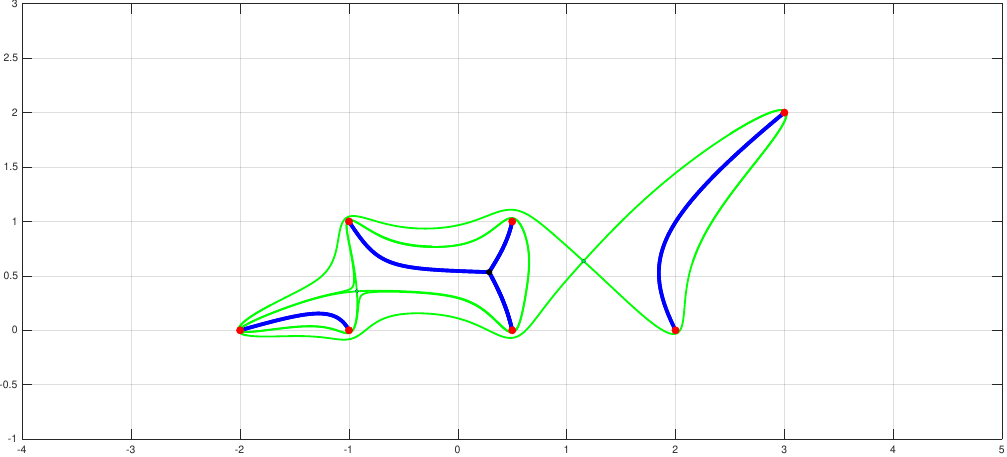} 
{\small Chebotarov continuum for the set $E=\{-1+0i,0.5+0i,0.5+1i,-1+1i,2+0i,-2+0i,3+2i\}$ with $L=2$.} 
\end{minipage} 
\caption{In these two examples of the generalized Chebotarov problem  where the set $\E$ is the same and we have two distinct solutions of the Boutroux problem \ref{problemPhiE} with $L=2$ stagnation points.}
\label{figunonuq}
\end{figure}

There is no implied statement about the {\it uniqueness} of the minimum. In fact it is easy to devise situations where there are several different absolute minima (all at $\scr F=0$). An example (see discussion of the generalized Chebotarov problem in section 6.2 of  \cite{MFRakh}) is for Chebotarov problems with $L>0$, where the minimal capacity set for a given set $\E$ is allowed to have  up to $L+1$ connected components; in this case there might be several solutions of the Boutroux problem. See Figure \ref{figunonuq}. In fact it is a possibly quite hard problem to compute the number of solutions; to understand the reason of the difficulty we need to keep in mind that $\scr F$ is not an analytic function and thus there are no degree considerations that are effective; even the  analogous  problem of estimating the number of real solutions for a system of polynomial equations is not possible to solve exactly. It is conceivable that one should be able to obtain an upper bound on the number of solutions, but this transcends the scope of this paper.

The convenient feature of this approach, however, is that from a numerical perspective it is simple to {\it verify} if we have reached a solution, since it suffices to compute the value of the functional $\scr F$ (or the individual integrals) and check whether they are numerically close to zero. 

\paragraph{Gradient descent analysis.}
We  denote with $\dot Q, \dot S, $ etc. the infinitesimal deformation of the quantities, noting that $\dot E(z)\equiv 0$. 
Let us denote by $ \gamma_1,\dots, \gamma_{2g}$ any basis of homology for the Riemann surface $\scr R$; in practice these are bounded closed contours containing an even number of branch-points of $\scr R$, typically just $2$. 

Let us choose a set of branch-cuts $\frak B$ for $\sqrt{Q}$ and denote the corresponding radical 
\be
\mathcal R(z) = \sqrt{Q(z)}: \C\setminus \frak B \to \C,
\ee
where the determination is fixed by the requirement that $\mathcal R(z) = t_R z^{R-1}  (1 + \mathcal O(z^{-1})$ as $z\to \infty$.

Let us denote 
\bea
T_\ell[Q]&:= \oint_{\gamma_0}  z^{-\ell} \sqrt{Q(z)} \frac{\d z}{2i\pi}, \ \ \ell = 0,1,\dots, R,\\
 P_j [Q]&:= \oint_{\gamma_j} \sqrt{Q(z)} \d z, \ \ j  = 1,2, \dots, 2g.
\eea
Here $\gamma_0$ is a large counterclockwise circle.
Then $\scr F$ writes
\be
\scr F[Q] = \frac 1 2 \sum_{\ell=0}^R |T_\ell[Q]-t_\ell|^2 + \frac 1 2 \sum_{\ell=1}^{2g} \Re(P_\ell[Q])^2,
\ee
and  the infinitesimal variation of $\scr F$ is 
\be
\frac {\d}{\d t}{\scr F[Q]} = 
\sum_{\ell=0}^R\Re \le[\ov{\le(T_\ell - t_\ell\ri)}\oint_{\gamma_0} \frac {z^{-\ell} \dot Q}{4i\pi \sqrt{Q}} \d z\ri] + \sum_{\ell=1}^{2g} \Re(P_\ell)
\Re\le(\oint_{\gamma_\ell} \frac {\dot Q  \d z}{2\sqrt{Q}}\ri).
\ee
Thus the goal of the gradient descent algorithm is to choose $\dot Q$ in such a way that 
$$
\frac {\d}{\d t}{\scr F[Q]} =-\scr F[Q].
$$
In components this means that $\dot Q$ should satisfy  
\bea
\label{flow0}
\oint_{\gamma_0} \frac {\dot Q}{4i\pi z^\ell \sqrt{Q}} \d z  &= t_\ell- T_\ell[Q],\ \ \ \ \ \ell =0,1,\dots, R\\
\Im\le(\oint_{\gamma_\ell} \frac {\dot Q  \d z}{2\sqrt{Q}}\ri)&=-\Im\le(\oint_{\gamma_\ell} \sqrt{Q} \d z\ri), \ \ \ell = 1,\dots, 2g.
\eea
Note that $T_R[Q] = t_R$ as per our choice of $Q$ in \eqref{formQ}, so that the corresponding equation in \eqref{flow0} is trivially satisfied. 
The deformation $\dot Q$ should be taken to preserve the form of $Q$, namely the number of stagnation points (zeros of $S$). The terminology refers to the ``orthogonal flow'', that is,  the stationary points of the foliation of vertical trajectories.

\paragraph{Initialization.}
  Choose $L, M$ with $2L+M = 2R+N-2$.  Choose a random $Q(z)$ of the form \eqref{formQ}. 
  
\paragraph{Main loop.}
Given the form \eqref{formQ} we must have 
\be
\dot {\sqrt{Q}} =\frac {\dot Q}{2\sqrt{Q}} = \frac{\dot S(z) \Delta(z) + \frac 1 2 S(z) \dot \Delta(z)}{\sqrt{\Delta(z) E(z)}} =: \frac {\dot F(z)}{\sqrt{\Delta(z) E(z)}}
\ee
We have $\deg(\dot S(z))\leq L-1$ since the leading coefficient of $S$ is fixed.  Thus the polynomial $\dot F(z)$ is of degree at most $\deg \dot F\leq M+L-1 = g+R-1 $.
We write it as follows
\be
\label{dotF}
\dot F(z) = \sum_{\ell=0}^{g+R-1} {f_\ell} z^\ell.
\ee
The coefficients $ {f_\ell}$ appear {\it linearly} in the  system \eqref{flow0}:
\bea
\label{sysf0}
 \sum_{k=0}^{g+R-1} {f_k} \oint_{\gamma_0} \frac { z^{k-\ell}\d z}{2i\pi\sqrt{\Delta(z) E(z)}}&= t_\ell-T_\ell[Q]\ \ \ \ell = 0,\dots, R-1\nn
 \\
  \sum_{k=0}^{g+R-1} \Re\le[ {f_k}  \oint_{\gamma_\ell} \frac { z^k\d z}{\sqrt{\Delta(z) E(z)} }\ri] &= -\Re P_\ell[Q],\ \ \ \ell =1,\dots,2g.
\eea

The system is determined; the reason is that  if there was a solution $[ {f_0} \dots,  { f_{g+R-1}}] $ of the homogeneous system then this would imply that the differential $\omega = \frac{\dot F \d z}{\sqrt{\Delta E}}$ has no pole at infinity on account of the first $R$ equations and is thus holomorphic. Then the remaining $2g$ equations state that all the periods are imaginary and a well--known consequence of the Riemann bilinear identities implies that there is no nontrivial holomorphic differential with purely imaginary periods. 

It is expedient to consider the system for the real and imaginary parts of the vector $\mathbf f = [f_0,\dots, f_{g+R-1}]$ so that the system becomes of size $(2g+2R)^2$.

Let us define the extended period matrix of size $(2g+R)\times(g+R)$
\be
\mathbb P_{a,b}:= \le\{
\begin{array}{cl}
\ds \oint_{\gamma_0} \frac { z^{b-a}\d z}{2i\pi \sqrt{\Delta(z) E(z)}} & {b= 1,\dots, g+R, \atop  a=1\dots, R}\\[10pt]
\ds \oint_{\gamma_{a-R}} \frac { z^{b-1} \d z}{\sqrt{\Delta(z) E(z)}} &{ b= 1,\dots, g+R, \atop  a=R+1\dots, 2g+R}
\end{array}
\ri.
\ee
and split it for convenience into the top $R\times (g+R)$ block, $\mathbb A$ and the bottom $2g\times (g+R)$ block, $\mathbb B$:
\be
\mathbb P = \le[\begin{array}{c}
\\
\mathbb A\in {\rm Mat}_{R\times (g+R)}(\C)\\[10pt]
\hline\\ \mathbb B\in {\rm Mat}_{2g\times (g+R)}(\C)
\\
\phantom{4}
\end{array}\ri].
\ee
Then the matrix of the system \eqref{sysf0} for the vector   $\mathbf x = [\Re \mathbf f, \Im \mathbf f]^t$ is 
\be
\le[
\begin{array}{c|c}
\Im \mathbb A & \Re \mathbb A\\
\hline
\Re \mathbb A & -\Im \mathbb A\\
\hline 
\Re \mathbb B & - \Im \mathbb B
\end{array}
\ri] \le[
\begin{array}{c}
\Re \mathbf f
\\
\Im \mathbf f
\end{array}\ri] = \le[
\begin{array}{c}
{}[\Im (t_\ell - T_\ell)]_{\ell=0}^{R-1}\\[10pt]
\hline
{}[\Re (t_\ell - T_\ell)]_{\ell=0}^{R-1}\\[10pt]
\hline 
{}[-\Re P_\ell]_{\ell=0}^{2g}
\end{array}
\ri].
\label{sysf00}
\ee
Thus the vector $\mathbf f = [f_0,\dots, f_{g+R-1}]^t$ is found by solving the regular linear system \eqref{sysf00}. 
At this point we have determined $\dot F$ and we need to determine $\dot S, \dot \Delta$ from 
\be
{\dot S(z) \Delta(z) + \frac 1 2 S(z) \dot \Delta(z)} =: {\dot F(z)}
\ee
Here $S,\Delta$ are known polynomials and, at this point so is $\dot F$. Once more, the system is a linear system for the coefficients of $\dot \Delta$ of degree $M-1$ and $\dot S$ of degree $L-1$. 

The matrix of this linear system is the Sylvester matrix of the polynomials $S,\frac 1 2\Delta$ and it is invertible as long as the two polynomials are relatively prime, which is our standing genericity assumption. 

Once $\dot S, \dot \Delta$ are determined we can implement a simple Euler integration method with small step $\d t$ (which can be adjusted heuristically as we progress in the evolution) 
\be
S \mapsto S +  \dot S\,\d t, \ \ \ \Delta \mapsto \Delta +\dot \Delta\, \d t .
\ee

\paragraph{Exit and other considerations.}
We continue to iterate the main loop until $\scr F$ has dropped below a pre-defined desirable threshold near $\scr F=0$ (typically $10^{-8}$, see also discussion in appendix). When $\scr F=0$  it means that the quadratic differential has the desired singular part and purely real periods. 

Of course the implementation of the algorithm could benefit from a more refined integration step, for example of Runge--Kutta  type, but we found that in practice it is sufficiently stable and fast. 

There are, of course, some potential numerical instabilities that may occur when the polynomials $S, \Delta$ share some root, but again we found rarely this is an issue. Depending on special situations, one can implement a ``merging check'' whereby if two roots of $\Delta$ are very close, we ``merge'' them and re-define $S$ to include the new simple root, $L\mapsto L+1$. 
\subsection{Comments on implementation}
We found that there are several small difficult details for a successful implementation, and the coding-oriented reader may benefit from the following pointers. We implemented the algorithm in Matlab and the commented code is available on GitHub \cite{github}. 

\paragraph{Choice of branch-cuts and definition of the radical.}
The computation of the various periods of $\sqrt{Q}$ requires to properly define  the function so that it has the branch--cuts where we expect them to be.  First of all we choose the $g+1$ branch-cuts as collection of pairwise non-intersecting straight segments, $\frak B$,  joining the $2g+2$ branchpoints (the zeros of $E$ and $\Delta$). This is effectively accomplished by the following recursive algorithm
\begin{enumerate}
\item {\bf Init:} Let $\scr Q_0$ be the set of all branch points and $\frak B=\emptyset$. Note that $\scr Q_0$ contains an even number of points $2g+2$;
\item  {\bf While $\scr Q_0\neq\emptyset$.} Compute the convex hull of $\scr Q_0$. Its boundary is a convex polygon. We add to $\frak B$ every second side of this polygon and remove the corresponding points from $\scr Q_0$. Repeat.
\end{enumerate}
The algorithm is not completely robust because if the convex hull is very ``thin'', the branch-cuts tend to almost overlap with resulting instability. Of course one may improve on this very basic approach to handle these situations. 

Once the branch-cuts have been defined, we need to define properly the radical $\mathcal R(z) = \sqrt {Q(z)}$. All numerical packages define the square root with a branch-cut on the negative axis. Thus if we want to define, for example, the radical $\sqrt{ \prod_{j=1}^{g+1} (z-a_j)(z-b_j)}$ with  straight cuts $[a_j,b_j]$, $j=1,\dots, g+1$, we need to code it as follows 
\be
\label{220}
\mathcal R(z) = \prod_{j=1}^{g+1} \frac{ \sqrt{ (\ov{a_j-b_j}) (z-a_j)}\sqrt{(\ov{a_j-b_j})(z-b_j)}}{\ov{(a_j- b_j)}}.
\ee
\paragraph{Integration over $\scr R$.}
The next annoying detail is that integration paths will inevitably cross the cuts, where our $\mathcal R(z)$ is discontinuous. Thus, coding a path $\gamma$ as a closed polygonal line, we need to check every segment of $\gamma$ and verify if it intersects any of the branchcuts $\frak B$. 
We consequently have to separate $\gamma$ into the parts on the two separate sheets of $\sqrt{Q}$ and integrate appropriately. This is necessary because no standard numerical package out of the box computes the integration along the analytic continuation.

\subsection{Special case: Chebotarov problems}
The Chebotarov problem consists in $t_0=1$ and $\Phi = 0$ ($R=0$) so that we have only the points of the set $\E$ and the choice of number of stagnation points.

\section{Pictures at an Exhibition}
We present here several outputs of the code. The green arcs are the critical trajectories from a stagnation point (zero of $S$). The red dots are points in $\E$, the black dots are zeros of $D$. 

\begin{minipage}{0.48\textwidth} 
\includegraphics[width=1\textwidth]{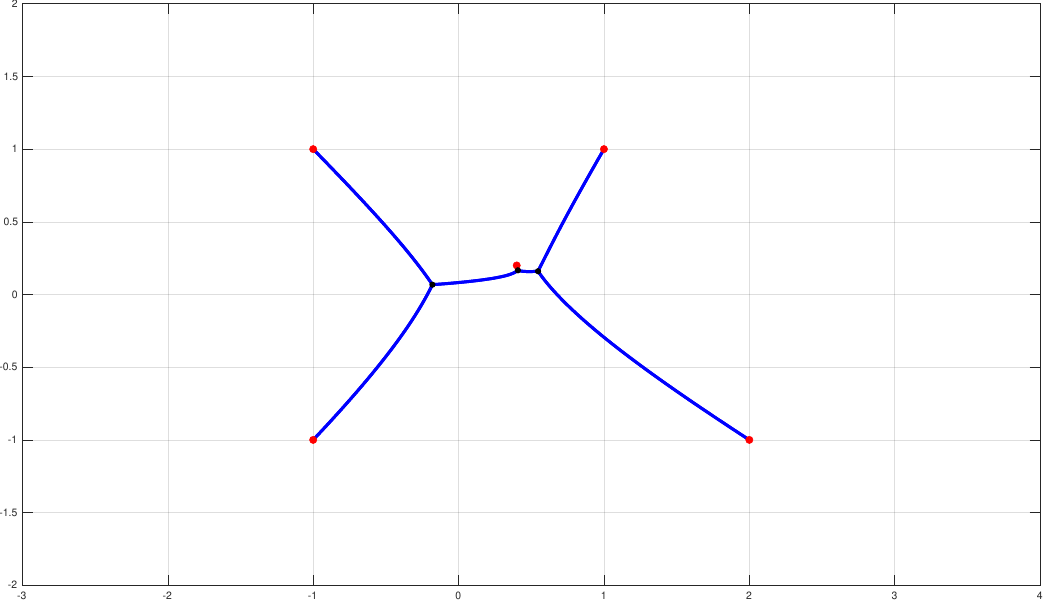} 
\captionof{figure}{\small Chebotarov continuum for the set $\E=\{-1+1i,-1-1i,0.4+0.2i,2-1i,1+1i\}$ with $L=0$.} 
\end{minipage} 
\hfill
\begin{minipage}{0.48\textwidth} 
\includegraphics[width=1\textwidth]{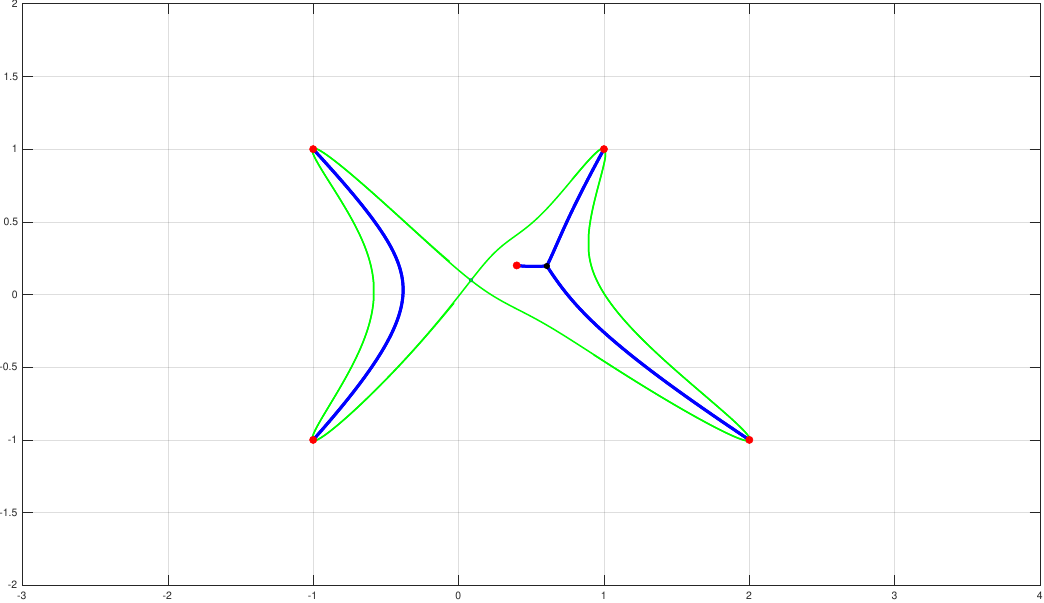} 
\captionof{figure}{\small Chebotarov continuum for the set $\E=\{-1+1i,-1-1i,0.4+0.2i,2-1i,1+1i\}$ with $L=1$.} 
\end{minipage}

\begin{minipage}{0.48\textwidth} 
\includegraphics[width=1\textwidth]{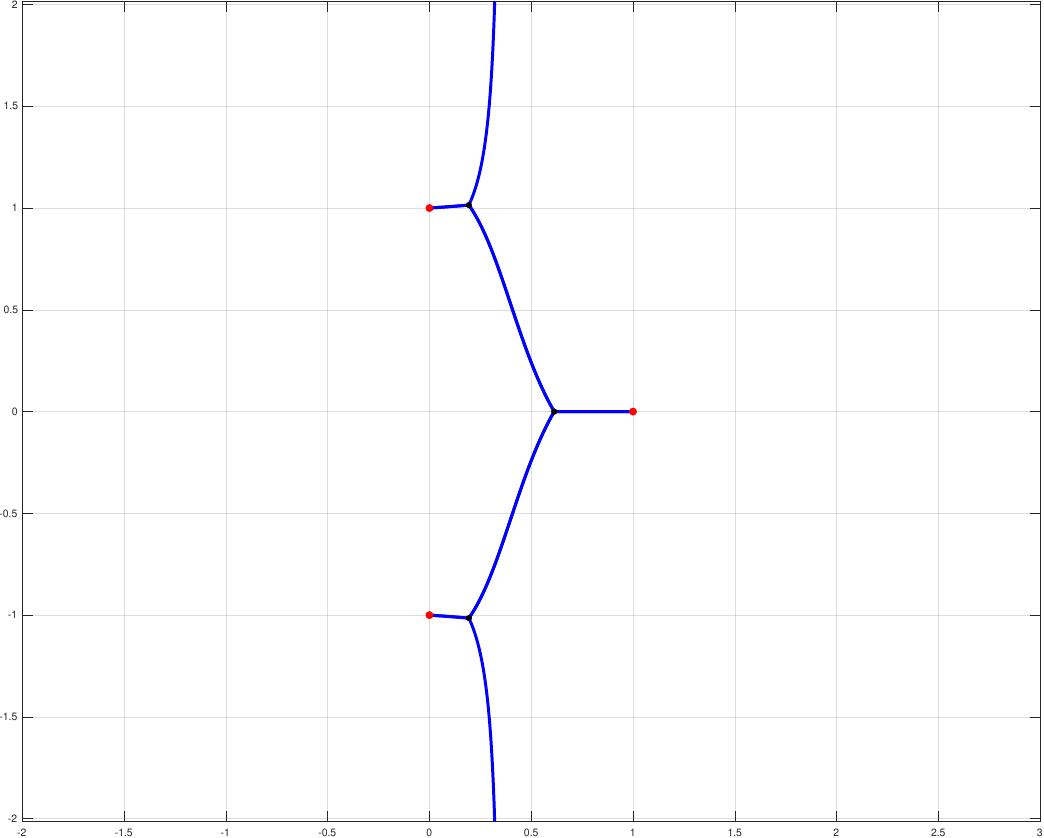} 
\captionof{figure}{\small In this case $\Phi = \sum_{\ell=1}^1 t_\ell z^\ell$ with $[t_1,\dots,t_1]=[1]$. Here  $t_0=0$. The set $\E=\{0-1i,0+1i,1+0i\}$ and $L=0$.} 
\end{minipage}\hfill
\begin{minipage}{0.48\textwidth} 
\includegraphics[width=1\textwidth]{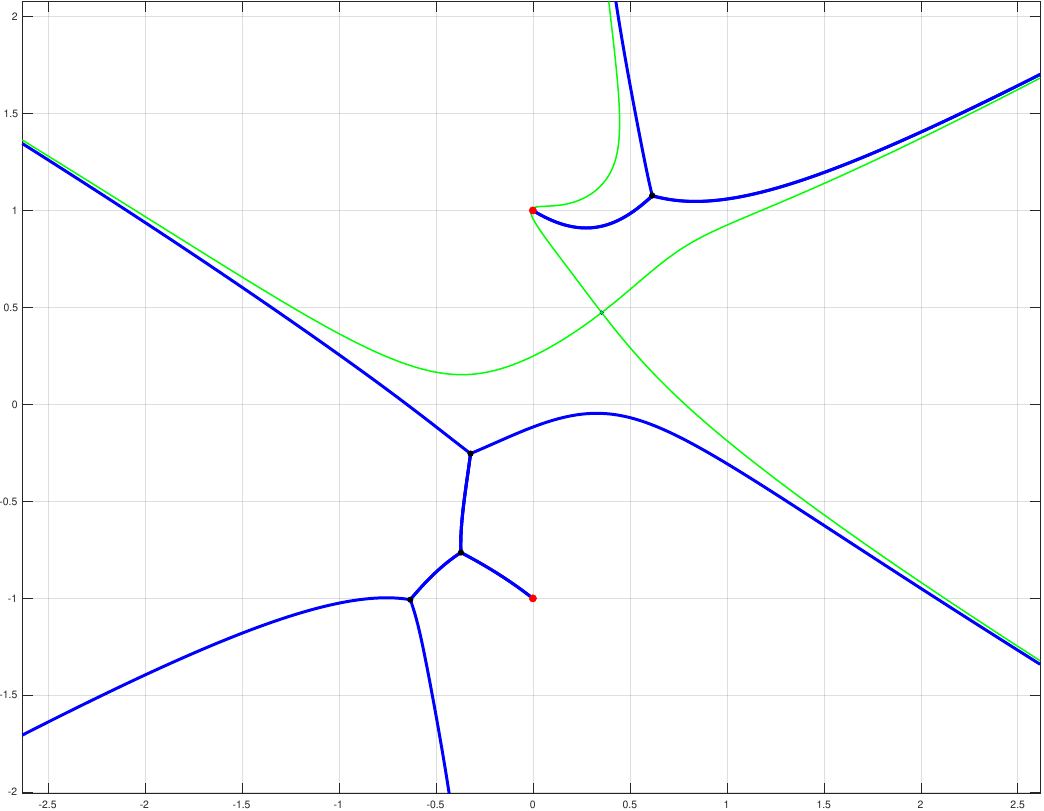} 
\captionof{figure}{\small In this case $\Phi = \sum_{\ell=1}^3 t_\ell z^\ell$ with $[t_3,\dots,t_1]=[1+0i,0+0i,0-1i]$. Here  $t_0=0$. The set $\E=\{0-1i,0+1i\}$ and $L=1$.} 
\end{minipage}

\begin{minipage}{0.45\textwidth} 
\includegraphics[width=1\textwidth]{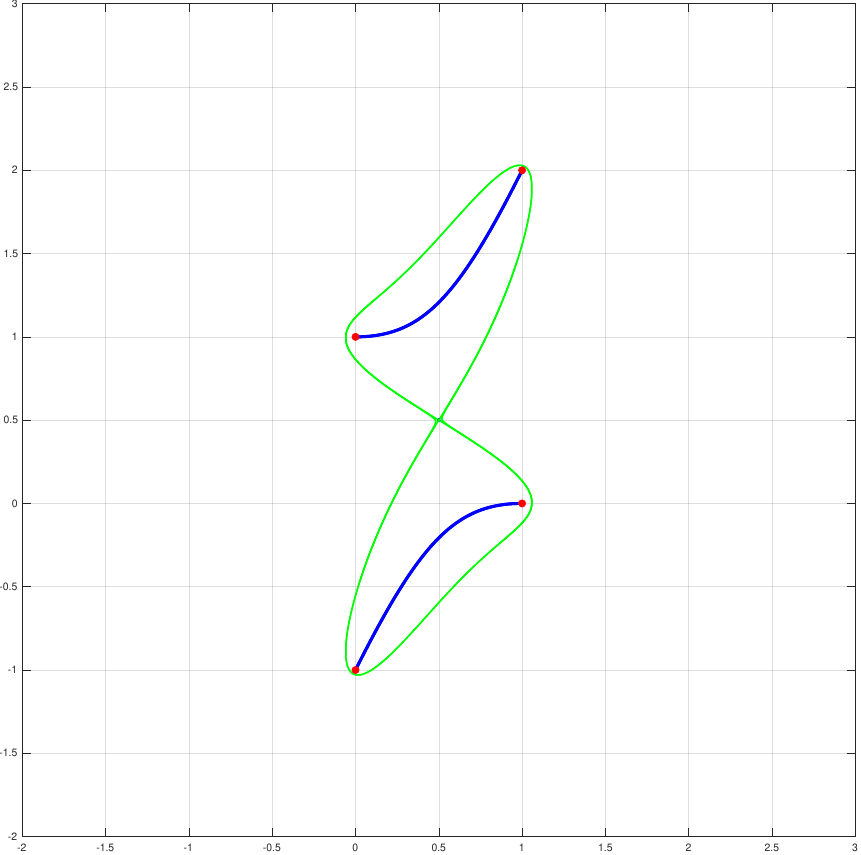} 
\captionof{figure}{\small Chebotarov continuum for the set $\E=\{0-1i,0+1i,1+0i,1+2i\}$ with $L=1$.} 
\end{minipage} \hfill
\begin{minipage}{0.48\textwidth} 
\includegraphics[width=1\textwidth]{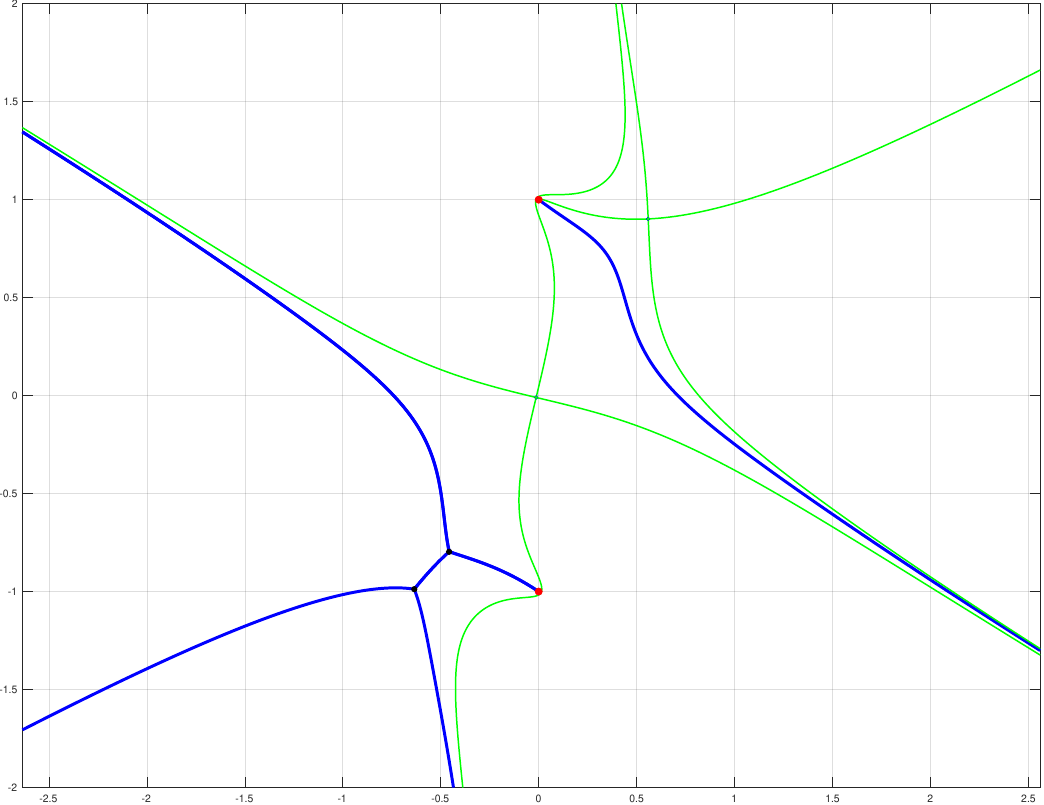} 
\captionof{figure}{\small In this case $\Phi = \sum_{\ell=1}^3 t_\ell z^\ell$ with $[t_3,\dots,t_1]=[1+0i,0+0i,0-1i]$. Here  $t_0=0$. The set $\E=\{0-1i,0+1i\}$ and $L=2$.} 
\end{minipage}

\begin{minipage}{0.45\textwidth} 
\includegraphics[width=1\textwidth]{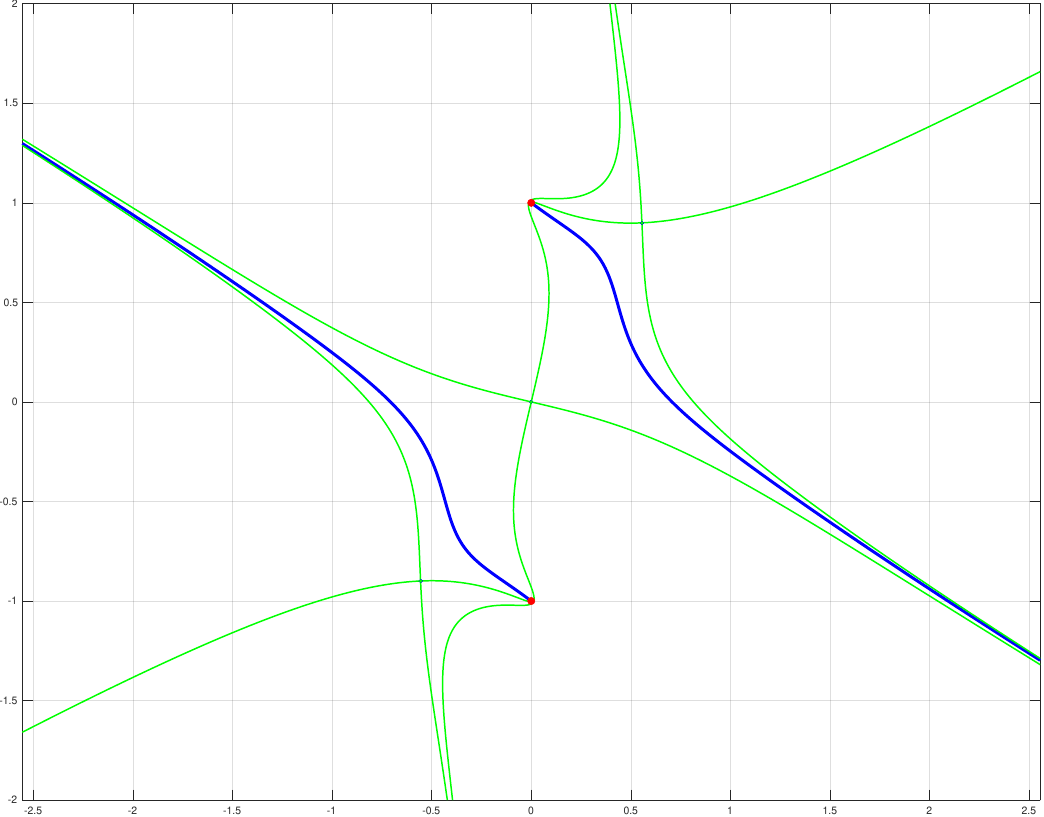} 
\captionof{figure}{\small In this case $\Phi = \sum_{\ell=1}^3 t_\ell z^\ell$ with $[t_3,\dots,t_1]=[1+0i,0+0i,0-1i]$. Here  $t_0=0$. The set $\E=\{0-1i,0+1i\}$ and $L=3$.} 
\end{minipage}\hfill
\begin{minipage}{0.45\textwidth} 
\includegraphics[width=1\textwidth]{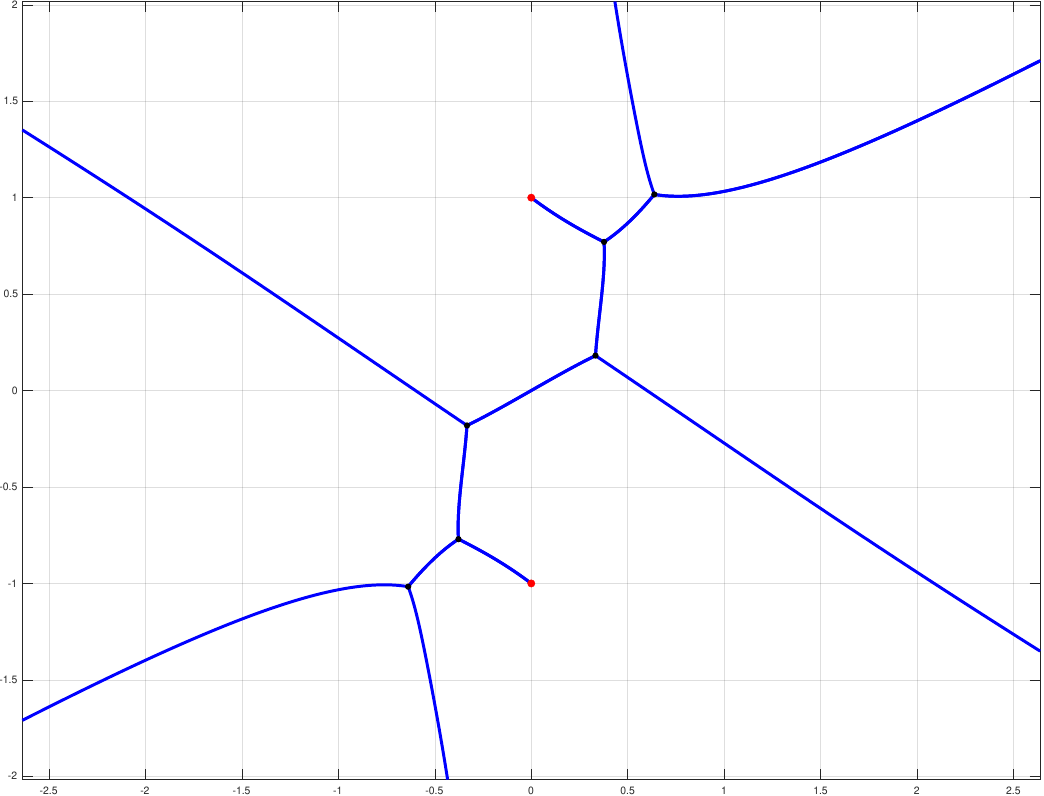} 
\captionof{figure}{\small In this case $\Phi = \sum_{\ell=1}^3 t_\ell z^\ell$ with $[t_3,\dots,t_1]=[1+0i,0+0i,0-1i]$. Here  $t_0=0$. The set $\E=\{0-1i,0+1i\}$ and $L=0$.} 
\end{minipage}

\begin{minipage}{0.45\textwidth} 
\includegraphics[width=1\textwidth]{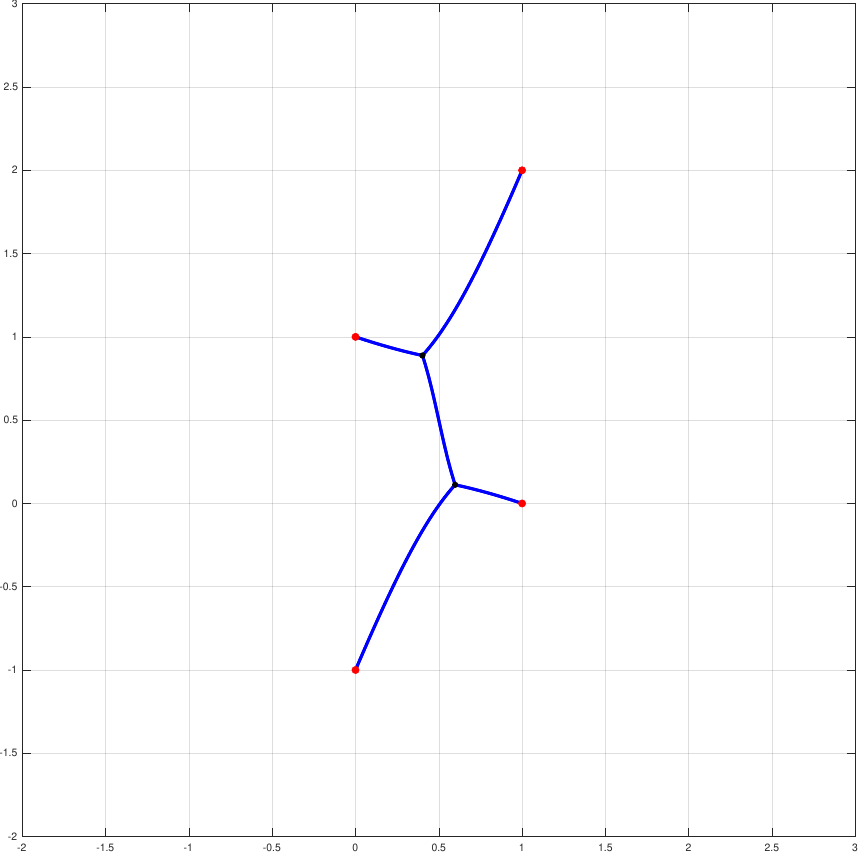} 
\captionof{figure}{\small Chebotarov continuum for the set $\E=\{0-1i,0+1i,1+0i,1+2i\}$ with $L=0$.} 
\end{minipage}\hfill
\begin{minipage}{0.5\textwidth} 
\includegraphics[width=1\textwidth]{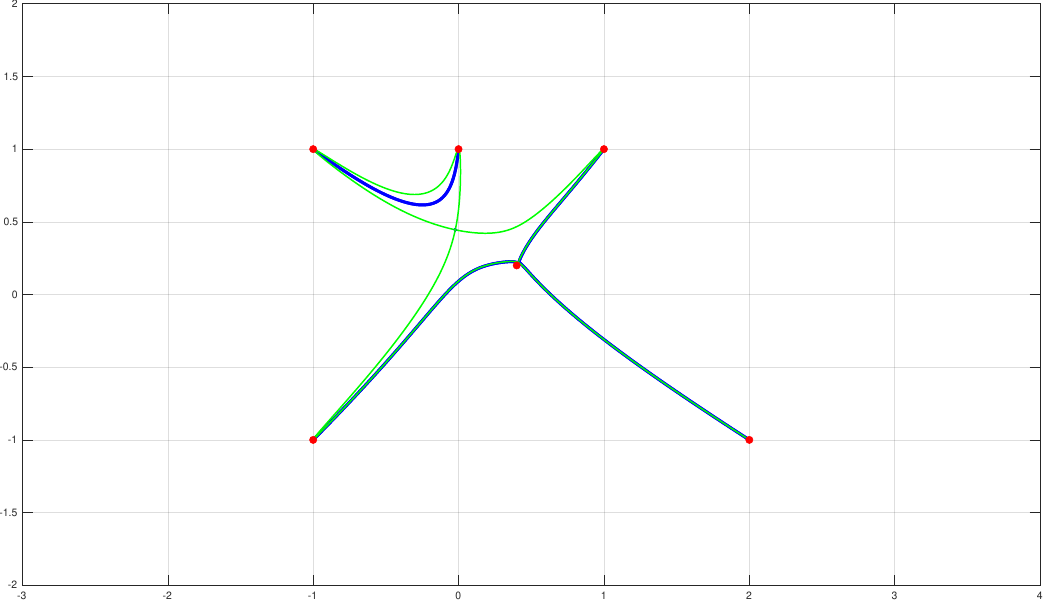} 
\captionof{figure}{\small Chebotarov continuum for the set $\E=\{-1+1i,-1-1i,0.4+0.2i,2-1i,1+1i,0+1i\}$ with $L=2$.} 
\end{minipage} 

\begin{minipage}{0.5\textwidth} 
\includegraphics[width=1\textwidth]{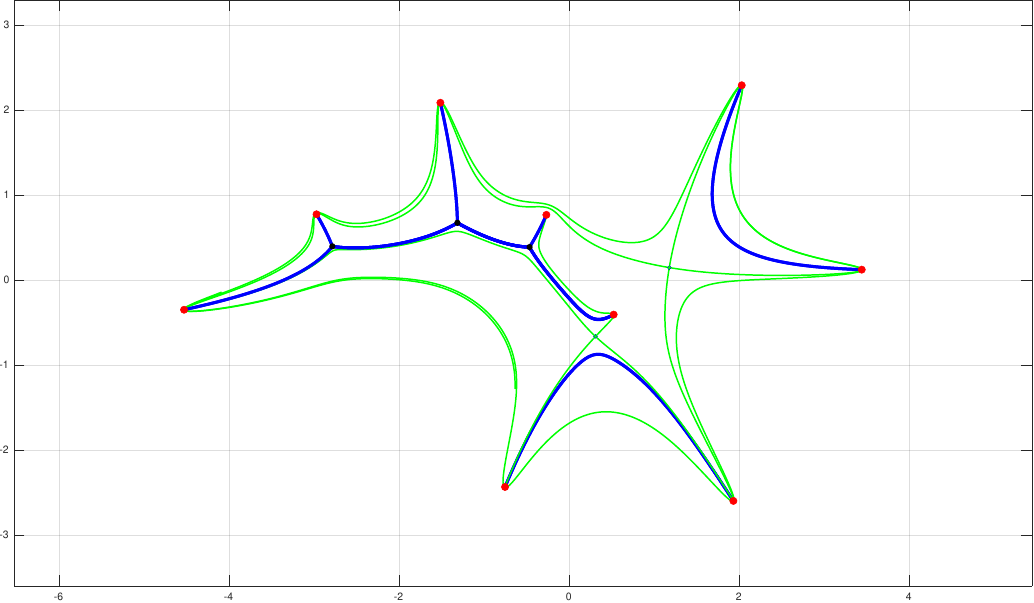} 
\captionof{figure}{\small Chebotarov continuum for the set $\E=\{3.4437+0.12426i,1.9331-2.5966i,2.0317+2.2937i,-4.5294-0.34736i,-0.75302-2.4321i,-2.9711+0.7759i,-1.5137+2.0879i,-0.26703+0.76894i,0.52499-0.40457i\}$ with $L=2$.} 
\end{minipage} \hfill
\begin{minipage}{0.5\textwidth} 
\includegraphics[width=1\textwidth]{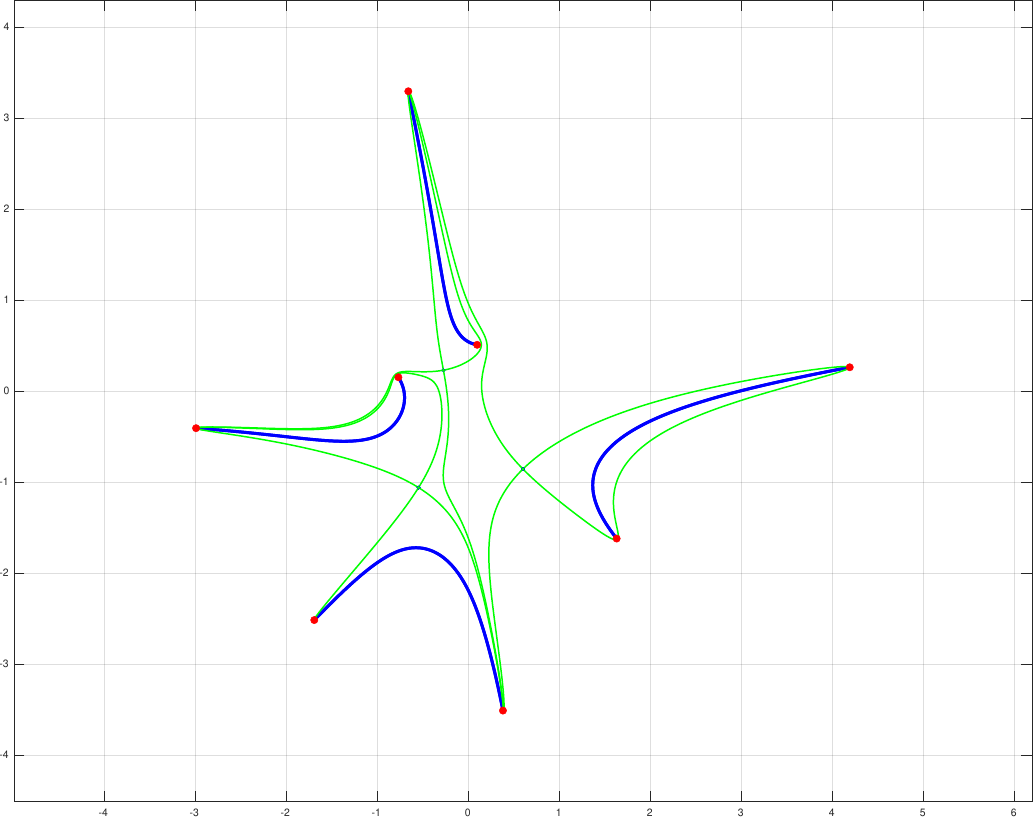} 
\captionof{figure}{\small Chebotarov continuum for the set $\E=\{4.1959+0.26274i,-0.65505+3.2958i,0.38502-3.5097i,1.6344-1.6195i,-2.9873-0.40647i,-1.6882-2.5151i,0.10121+0.50947i,-0.76412+0.1532i\}$ with $L=3$.} 
\end{minipage} 
\appendix
\section{The code package}
The algorithm was implemented in Matlab and the package consists of the files below, available on Github \cite{github}. 
We briefly comment on their usage, but the input/output arguments may be changed easily by editing the files (the description below may differ from future revisions).
\begin{enumerate}
\item ``G\_functions.m''; the main arguments  are  $\Phi$ and $\E$ as in Problem \ref{problemPhiE}. Additionally one can specify $L$, the degree of $S$, i.e. the initial desired  number of stagnation points. The main output consists of the polynomials $S$ and $\Delta$.  Polynomials are specified in the standard Matlab way, namely, as list of coefficients. Note however that the constant coefficient passed via $\Phi$ is interpreted as $t_0$ (the residue): see  Problem \ref{problemPhiE}.  For $\Phi=0$ it produces the solutions of the Chebotarov problems, the original version if $L=0$ and the generalized poly-continua if $L>0$. Specific details are in the help to the function itself.
\item ``CriticalTrajectories\_Rational.m''; this is a function used to plot the critical trajectories of a rational quadratic differential. It works with any quadratic differential of the form $\mathcal Q(z)= Q(z) \d z^2$, with $Q(z)$ an arbitrary rational function and will plot all the critical horizontal trajectories. The computation of each trajectory is done by ``QTrajectories.m'' below.  See function help for usage. 
\item ``QTrajectories.m''; it computes the critical trajectory of a rational quadratic differential using a Runge-Kutta 4 method. It is invoked by (2.) for the plotting.
\item ``Hyperelliptic\_Integral.m''; this function computes the integral of a radical function $F(z)$ along a specified path (given as a sequence of point in the plane specifying a polygonal line). The argument are the handle to a function object $F$, the path and a list of cuts (specified as pair of points) that indicate the endpoints of straight segments where the cuts of the radical $F$ are placed. No attempt is made to verify that the function and the cuts are consistent with each other. The radical should have been defined as explained around \eqref{220} and the cuts should be passed consistently with that definition. The output is the value of the integral of the analytic continuation of $F$ along the given path, keeping into account that the path may jump between the two sheets of $F$. The function is invoked by (1.). More details are in the comments inside the file.
\end{enumerate}    
The main file is $(1.)$ and it calls the others.

\subsection{Error and performance control}
The author of the present paper is at best describable as an amateur coder and the focus of the paper is on the algorithm rather than the implementation. Since we are finding numerically the zeros of $\scr F$ defined by \eqref{defF}, the upper bound on absolute values of the single integrals is given simply by the square root of $\scr F$. The precision in the evaluation of the integrals is the one provided natively by the numerical implementation of the Matlab command ``integral'',
which we used with the default absolute error tolerance (AbsTol$=10^{-10}$) and relative error tolerance (RelTol$=10^{-6}$).
In the code we actually compute $Fun = \sqrt{\scr F}$  and given the default absolute error tolerance we can set a meaningful threshold as $10^{-9}$ by default. 

The implementation is quite fast; for all the examples produced below it takes less than a few seconds and less than hundred iterations. Of course there are variations as the most expensive part of the iteration is the computation of the several integrals. 
Table \ref{table}  gives some indicative values of the performance of a typical  run to achieve the default threshold $\sqrt{\scr F}\leq 10^{-9}$. Since the initial conditions are chosen at random, each run has a fluctuation in run-time even if started with the same parameters. We indicate the result of four runs and corresponding times for $0\leq L\leq N/2-1$ (the maximum number of stagnation points) in Table \ref{table}. 
\begin{table}
\begin{center}\begin{tabular}{|c|c|c|c|c|c|}
\hline
$N$ & $g$ &$L$ & ITER & Time(sec)\\
\hline
18 & 16 & 0 & \{32, 31, 32, 31\} & \{49.5, 54.1, 61.8,  62.5\}\\
\hline 
18 & 15 & 1 & \{37, 37, 40, 41\}& \{51.6, 54.0, 60.0, 65.3\}\\
\hline
18 & 14 & 2 & \{38, 39, 46, 46\}& \{50.5, 52.2, 54.4, 60.0\}\\
\hline
18 & 13 & 3 & \{40, 38, 46, 38\}& \{42.8, 42.7, 46.8, 44.1\}\\
\hline
18 & 12 & 4 & \{40, 42, 45, 41\}& \{40.1,  44.1, 41.2, 42.9\}\\
\hline
18 & 11 & 5 & \{47, 45, 43, 47\}& \{35.7, 38.6, 37.1, 36.7\}\\
\hline
18 & 10 & 6 & \{42, 43, 44, 41\}& \{28.5, 30.5, 32.2, 29.9\}\\
\hline
18 & 9 & 7 & \{35,  41, 42, 43\} &\{19.1,  24.0,  26.9, 26.2\}\\
\hline
18 & 8 & 8 & \{33, 33, 32,  32\}& \{15.9, 16.9, 16.0, 16.5\}\\
\hline
\end{tabular}\end{center}
\caption{Some indicative values of the number of iterations (ITER) and  the total time (including the plotting) of a typical  run to achieve the default threshold $\sqrt{\scr F}\leq 10^{-9}$ (run on a commodity desktop machine). The column $g$ indicates the genus.}
\label{table}
\end{table}

\end{document}